\documentclass[]{amsproc}

\usepackage{aspmproc}         
\usepackage{amssymb}

\usepackage{graphicx,psfrag}  
 
\usepackage[dvipdfm,bookmarks=true,bookmarksnumbered=
true,bookmarkstype=toc]{hyperref}   

\rcvdate{January 15, 2009}
\rvsdate{March 3, 2009}   

\numberwithin{equation}{section}
\numberwithin{figure}{section}

 \newtheorem{thm}{Theorem}\numberwithin{thm}{section}
  \newtheorem{lem}{Lemma}\numberwithin{lem}{section}
  \numberwithin{slem}{section} 
   \numberwithin{cor}{section}
  \numberwithin{prop}{section}
  \numberwithin{fact}{section}

\theoremstyle{definition}
  \newtheorem{defn}{Definition}\numberwithin{defn}{section}
  \numberwithin{ex}{section}

\theoremstyle{remark}
\newtheorem{rem}{Remark}\numberwithin{rem}{section}

\newtheorem*{ack}{\bf Acknowledgment}

\numberwithin{equation}{section}

\DeclareMathOperator{\vol}{vol}
\DeclareMathOperator{\Ric}{Ric}
\DeclareMathOperator{\BG}{BG}
\DeclareMathOperator{\CD}{CD}

\DeclareMathOperator{\ap}{ap}

\newcommand{\Hm}{\mathcal{H}}
\newcommand{\Cut}{\text{\rm Cut}}

\newcommand{\field}[1]{\mathbb{#1}}

\newcommand{\R}{\field{R}}

\title[Bishop-Gromov condition]{Infinitesimal Bishop-Gromov condition\\
  for Alexandrov spaces}

\author[K. Kuwae]{Kazuhiro Kuwae}
\address{Department of Mathematics and Engineering\\
  Graduate School of Science and Technology\\
  Kumamoto University\\ 
  Kumamoto, 860-8555, JAPAN}
\email{kuwae@gpo.kumamoto-u.ac.jp}

\author[T. Shioya]{Takashi Shioya}
\address{Mathematical Institute\\
 Tohoku University\\
  Sendai 980-8578, JAPAN}
\email{shioya@math.tohoku.ac.jp}

\rcvdate{January 15, 2009}
\rvsdate{2009}

\thanks{The authors are partially supported by a Grant-in-Aid
  for Scientific Research No.~19540220 and 20540058 from
  the Japan Society for the Promotion of Science}


\begin{document}

\begin{abstract}
  We prove the infinitesimal version of Bishop-Gromov volume comparison
  condition for Alexandrov spaces.
\end{abstract}

\maketitle

\section{Introduction} \label{sec:intro}

We first present the definition of the infinitesimal Bishop-Gromov
volume comparison condition for Alexandrov spaces.

For a real number $\kappa$, we set
\[
s_\kappa(r) :=
\begin{cases}
  \sin(\sqrt{\kappa}r)/\sqrt{\kappa}  &\text{if $\kappa > 0$},\\
  r &\text{if $\kappa = 0$},\\
  \sinh(\sqrt{|\kappa|}r)/\sqrt{|\kappa|} &\text{if $\kappa < 0$}.
\end{cases}
\]
The function $s_\kappa$ is the solution of the Jacobi equation
$s_\kappa''(r) + \kappa s_\kappa(r) = 0$ with initial condition
$s_\kappa(0) = 0$, $s_\kappa'(0) = 1$.

Let $M$ be an Alexandrov space and set $r_p(x) := d(p,x)$ for $p,x \in M$,
where $d$ is the distance function.
For $p \in M$ and $0 < t \le 1$,
we define a subset $W_{p,t} \subset M$ and a map $\Phi_{p,t} : W_{p,t} \to M$
as follows.
We first set $\Phi_{p,t}(p) := p \in W_{p,t}$.
A point $x$ ($\neq p$) belongs to $W_{p,t}$
if and only if there exists $y \in M$ such that
$x \in py$ and $r_p(x) : r_p(y) = t:1$, where
$py$ is a minimal geodesic from $p$ to $y$.
Since a geodesic does not branch on an Alexandrov space,
for a given point $x \in W_{p,t}$ such a point $y$ is unique and we set
$\Phi_{p,t}(x) := y$.
The triangle comparison condition
implies the local Lipschitz continuity of the map $\Phi_{p,t} : W_{p,t} \to M$.
We call $\Phi_{p,t}$ the \emph{radial expansion map}.

Let $\mu$ be a positive Radon measure with full support in $M$,
and $n \ge 1$ a real number.

\medskip
{\bf Infinitesimal Bishop-Gromov Condition $\BG(\kappa,n)$ for $\mu$:}\\
For any $p \in M$ and $t \in (\,0,1\,]$, we have
\[
d(\Phi_{p,t\,*} \mu)(x)
\ge \frac{t\,s_\kappa(t\,r_p(x))^{n-1}}{s_\kappa(r_p(x))^{n-1}}
d\mu(x)
\]
for any $x \in M$ such that $r_p(x) < \pi/\sqrt{\kappa}$ if $\kappa > 0$,
where $\Phi_{p,t\,*}\mu$ is the push-forward
by $\Phi_{p,t}$ of $\mu$.

\medskip
For an $n$-dimensional complete
Riemannian manifold, the Riemannian volume measure satisfies $\BG(\kappa,n)$
if and only if the Ricci curvature satisfies $\Ric \ge (n-1)\kappa$
(see Theorem 3.2 of \cite{Ota:mcp} for the `only if' part).
We see some studies on similar (or same) conditions to $\BG(\kappa,n)$ in
\cite{CC:strRicI,St:heat,KwSy:geneMCP,KwSy:sobmet,Rm:Poincare,Ota:mcp,
Wt:localcut} etc.
$\BG(\kappa,n)$ is sometimes called the Measure Contraction Property
and is weaker than the curvature-dimension (or lower $n$-Ricci curvature)
condition, $\CD((n-1)\kappa,n)$, introduced
by Sturm \cite{St:geomI,St:geomII} and Lott-Villani \cite{LV:Ricmm}
in terms of mass transportation.
For a measure on an Alexandrov space, $\BG(\kappa,n)$ is equivalent to
the $((n-1)\kappa,n)$-MCP introduced by Ohta \cite{Ota:mcp}.
In our paper \cite{KwSy:splitting,KwSy:lapcomp},
we prove a splitting theorem under $\BG(0,N)$.
For a survey of geometric analysis on Alexandrov spaces,
we refer to \cite{Sy:Ronsetsu}

The purpose of this paper is to prove the following

\begin{thm} \label{thm:BG}
  Let $M$ be an $n$-dimensional Alexandrov space of curvature $\ge \kappa$.
  Then, the $n$-dimensional Hausdorff measure $\Hm^n$ on $M$
  satisfies the infinitesimal Bishop-Gromov condition $\BG(\kappa,n)$.
\end{thm}

Note that we claimed this theorem in Lemma 6.1 of \cite{KwSy:geneMCP},
but the proof in \cite{KwSy:geneMCP} is insufficient.
The theorem also completes the proof of
Proposition 2.8 of \cite{Ota:mcp}.

For the proof of the theorem, we have the delicate problem that
the topological boundary of the domain $W_{p,t}$
of the radial expansion $\Phi_{p,t}$ is not necessarily of $\Hm^n$-measure zero.
In fact, we have an example of an Alexandrov space such that the
cut-locus at a point is dense (see Remark \ref{rem:cutlocus}),
in which case the boundary of $W_{p,t}$
has positive $\Hm^n$-measure.
This never happens for Riemannian manifolds.
To solve this problem,
we need some delicate discussion using
the approximate differential of $\Phi_{p,t}$.

\begin{ack}
  The authors would like to thank Professor Shin-ichi Ohta for
  his comments.
\end{ack}

\section{Preliminaries} \label{sec:prelim}

\subsection{Alexandrov spaces}

In this paper, we mean by an Alexandrov space a complete locally compact
geodesic space of curvature bounded below locally and
of finite Hausdorff dimension.
We refer to \cite{BGP,OS:rstralex,KMS:lap} for the basics
for the geometry and analysis on Alexandrov spaces.
Let $M$ be an Alexandrov space of Hausdorff dimension $n$.
Then, $n$ coincides with the covering dimension of $M$ which is a
nonnegative integer.
Take any point $p \in M$ and fix it.
Denote by $\Sigma_pM$ the space of directions at $p$,
and by $K_pM$ the tangent cone at $p$.
$\Sigma_pM$ is an $(n-1)$-dimensional compact Alexandrov space
of curvature $\ge 1$ and
$K_pM$ an $n$-dimensional Alexandrov space of curvature $\ge 0$.

\begin{defn}[Singular Point, $\delta$-Singular Point]
  A point $p \in M$ is called a \emph{singular point of $M$}
  if $\Sigma_pM$ is not isometric to the unit sphere $S^{n-1}$.
  For $\delta > 0$, we say that a point $p \in M$ is \emph{$\delta$-singular}
  if $\Hm^{n-1}(\Sigma_pM) \le \vol(S^{n-1})-\delta$.
  Let us denote the set of singular points of $M$ by
  $S_M$ and the set of $\delta$-singular points of $M$ by
  $S_\delta$.  
\end{defn}

We have $S_M = \bigcup_{\delta > 0} S_\delta$.
Since the map $M \ni p \mapsto \Hm^n(\Sigma_pM)$ is lower semi-continuous,
the set $S_\delta$ of $\delta$-singular points in $M$ is a closed set.

\begin{lem}[\cite{Pt:parallel}] \label{lem:geod}
  Let $\gamma$ be a minimal geodesic joining two points $p$ and $q$ in $M$.
  Then, the space of directions, $\Sigma_xM$, at all interior points of
  $\gamma$, $x \in \gamma \setminus \{p,q\}$, are isometric to each other.
  In particular, any minimal geodesic joining two non-singular
  {\rm(}resp.~non-$\delta$-singular{\rm)} points is contained in the set of
  non-singular {\rm(}resp.~non-$\delta$-singular{\rm)} points
  {\rm(}for any $\delta > 0${\rm)}.
\end{lem}

The following shows the existence of differentiable and
Riemannian structure on $M$.

\begin{thm}
  \label{thm:str}
  For an $n$-dimensional Alexandrov space $M$, we have the
  following:
  \begin{enumerate}
  \item  There exists a number $\delta_n > 0$ depending only on $n$
    such that
    $M^* := M \setminus S_{\delta_n}$ is a manifold
    {\rm(\cite{BGP})}
    and has a natural $C^\infty$ differentiable structure
    {\rm(\cite{KMS:lap})}.
  \item The Hausdorff dimension of $S_M$ is $\le n-1$
    {\rm(\cite{BGP, OS:rstralex})}.
  \item We have a unique continuous Riemannian metric $g$ on
    $M \setminus S_M \subset M^*$
    such that the distance function induced from $g$ coincides with
    the original one of $M$ {\rm(\cite{OS:rstralex})}.
    The tangent space at each point in $M \setminus S_M$ is
    isometrically identified with the tangent
    cone {\rm(}\cite{OS:rstralex}{\rm)}.
    The volume measure on $M^*$ induced from $g$
    coincides with the $n$-dimensional Hausdorff measure $\Hm^n$
    {\rm(\cite{OS:rstralex})}.
  \end{enumerate}
\end{thm}

\begin{rem}
  In \cite{KMS:lap} we construct a $C^\infty$ structure only on
  $M \setminus B(S_{\delta_n},\epsilon)$, where $B(A,\epsilon)$
  denotes the $\epsilon$-neighborhood of $A$.
  However this is independent of $\epsilon$ and extends to
  $M^*$.
  The $C^\infty$ structure is a refinement of the structures
  of \cite{OS:rstralex,Ot:secdiff,Pr:DC} and is
  compatible with the DC structure of \cite{Pr:DC}.
\end{rem}

Note that the metric $g$ is defined only on $M^* \setminus S_M$ and
does not continuously extend to any other point of $M$.

\begin{defn}[Cut-locus]
  Let $p \in M$ be a point.
  We say that a point $x \in M$ is a \emph{cut point of $p$}
  if no minimal geodesic from $p$ contains $x$ as an interior point.
  Here we agree that $p$ is a cut point of $p$.
  The set of cut points of $p$ is called the \emph{cut-locus of $p$} and
  denoted by $\Cut_p$.
\end{defn}

Note that $\Cut_p$ is not necessarily a closed set.
For the $W_{p,t}$ defined in \S\ref{sec:intro}, it follows that
$\bigcup_{0 < t < 1} W_{p,t} = X \setminus \Cut_p$.
The cut-locus $\Cut_p$ is a Borel subset and satisfies
$\Hm^n(\Cut_p) = 0$ (Proposition 3.1 of \cite{OS:rstralex}).

\begin{rem} \label{rem:cutlocus}
  There is an example of a $2$-dimensional Alexandrov space $M$ such that
  $S_M$ is dense in $M$ (see \cite{OS:rstralex}).
  For such an example, $\Cut_p$ for any $p \in M$ is also dense
  in $M$.
\end{rem}

\subsection{Approximate differential}

\begin{defn}[Density; cf.~2.9.12 in \cite{Fd:gmeas}]
  Let $X$ be a metric space with a Borel measure $\mu$.
  A subset $A \subset X$ \emph{has density zero at a point $x \in X$}
  if
  \[
  \lim_{r \to 0} \frac{\mu(B(x,r) \cap A)}{\mu(B(x,r))} = 0.
  \]
\end{defn}

\begin{defn}[Approximate Differential; cf.~3.1.2 in \cite{Fd:gmeas}]
  Let $A \subset \R^m$ be a subset and $f : A \to \R^n$ a map.
  A linear map $L : \R^m \to \R^n$ is called the
  \emph{approximate differential of $f$ at a point $x \in A$}
  if the approximate limit of
  \[
  \frac{|\,f(y)-f(x)-L(y-x)\,|}{|y-x|}
  \]
  is equal to zero as $y \to x$, i.e.,
  for any $\delta > 0$, the set
  \[
  \left\{\; y \in A \setminus \{x\}\; \Bigl|\;  
    \frac{|\,f(y)-f(x)-L(y-x)\,|}{|y-x|} \ge \delta\;\right\}
  \]
  has density zero at $x$, where we consider the Lebesgue (or equivalently
  $m$-dimensional Hausdorff) measure on $\R^m$ to measure the density.
  We say that \emph{$f$ is approximately differentiable at a point $x \in A$}
  if the approximate differential of $f$ at $x$ exists.
  Denote by `$\ap df_x$' the approximate differential of $f$ at $x$.
  It is unique at each approximate differentiable point.

  Let $M$ and $N$ be two differentiable manifolds and let $A \subset M$.
  We give a map $f : A \to N$ and a point $x \in A$.
  Take two charts $(U,\varphi)$ and $(V,\psi)$ around $x$ and
  $f(x)$ respectively.
  The map $f$ is said to be \emph{approximately differentiable at $x$}
  if $\psi\circ f\circ\varphi^{-1}$ is approximately differentiable
  at $\varphi(x)$.
  If $f$ is approximately differentiable at $x$, then
  the \emph{approximate differential} `$\ap df_x$' of $f$ at $x$
  is defined by 
  \[
  \ap df_x := (d\psi_{f(x)})^{-1}\circ\ap d(\psi\circ f\circ\varphi^{-1})_{\varphi(x)}
  \circ d\varphi_x : T_xM \to T_{f(x)}N.
  \]
  The approximate differentiability of $f$ at $x$
  and $\ap df_x$ are both independent of $(U,\varphi)$ and $(V,\psi)$.
\end{defn}

\section{Proof of Theorem \ref{thm:BG}}

Let $M$ be an Alexandrov space of curvature $\ge \kappa$.
We first investigate the exponential map on $M$.
Denote by $o_p$ the vertex of the tangent cone $K_pM$ at a point $p \in M$.
We denote by $U_p \subset K_pM$ the \emph{inside of the
tangential cut-locus of $p$}, i.e.,
$v \in U_p$ if and only if there is a minimal geodesic
$\gamma : [\,0,a\,] \to M$ from $p$ with $a > 1$
such that $\gamma'(0) = v$, where $\gamma'(t)$ denotes the element
of $K_{\gamma(t)}M$ tangent to $\gamma|_{[\,t,t+\epsilon\,)}$, $\epsilon > 0$,
and whose distance from $o_{\gamma(t)} \in K_{\gamma(t)}M$
is equal to the speed of parameter of $\gamma$.
Note that $U_p$ is not necessarily an open set.
Since the exponential map $\exp_p|_{U_p} : U_p \to M \setminus \Cut_p$
is a homeomorphism
and since $W_{p,t} \cap \bar B(p,r)$ is compact for any $0 < t \le 1$ and
$r > 0$, the set
\[
U_p = \bigcup_{0 < t \le 1,\ r > 0}(\exp_p|_{U_p})^{-1}(W_{p,t} \cap \bar B(p,r))
\]
is a Borel subset of $K_pM$.

Denote by $\Theta(t|a,b,\dots)$ a function of $t,a,b,\dots$
such that $\Theta(t|a,b,\dots) \to 0$ as $t \to 0$ for any fixed $a,b,\dots$.
We use $\Theta(t|a,b,\dots)$ as Landau symbols.

\begin{lem} \label{lem:HmKp}
  For any $p \in M$, $r > 0$, and
  for any $\Hm^n$-measurable subset $A \subset B(o_p,r) \subset K_pM$,
  we have
  \begin{align}
    |\,\Hm^n(\exp_p(A \cap U_p)) - \Hm^n(A)\,| &\le \Theta(r|p,n)\, r^n,
    \tag{1}\\
    \Hm^n(B(o_p,r) \setminus U_p) &\le \Theta(r|p,n)\, r^n.
    \tag{2}
  \end{align}
\end{lem}

Note that $\Theta(r|p,n)$ here is independent of $A$.

\begin{proof}
  Let $p \in M$ and $r > 0$.
  By the triangle comparison condition,
  $\exp_p : U_p \cap B(o_p,r) \to M$
  is Lipschitz continuous with Lipschitz constant $1+\Theta(r|p)$.
  Therefore, for any $\Hm^n$-measurable $A \subset B(o_p,r)$,
  \begin{align*}
    \Hm^n(A) &\ge (1-\Theta(r|p,n))\,\Hm^n(\exp_p(A \cap U_p)),\\
    \Hm^n(B(o_p,r) \setminus A) &\ge
    (1-\Theta(r|p,n))\,\Hm^n(B(p,r) \setminus \exp_p(A \cap U_p)).
  \end{align*}
  According to Lemma 3.2 of \cite{Sy:raysalex}, we have
  \[
  \lim_{\rho\to 0} \frac{\Hm^n(B(p,\rho))}{\rho^n} = \Hm^n(B(o_p,1))
  = \frac{\Hm^n(B(o_p,r))}{r^n}.
  \]
  Combining those three formulas we have the lemma.
\end{proof}

Let $p \in M$ and $0 < t \le 1$.
We restrict the domain of
the radial expansion map $\Phi_{p,t} : W_{p,t} \to M$
to the subset
\[
W'_{p,t} := W_{p,t} \setminus (\Phi_{p,t}^{-1}(\Cut_p) \cup S_{\delta_n}),
\]
where $S_{\delta_n}$ is as in Theorem \ref{thm:str}.

\begin{lem} \label{lem:Wprime}
  We have $\Phi_{p,t}(W'_{p,t}) = M \setminus (\Cut_p \cup S_{\delta_n})$
  and the map $\Phi_{p,t}|_{W'_{p,t}} : W'_{p,t} \to
  M \setminus (\Cut_p \cup S_{\delta_n})$
  is bijective.
  In particular, the sets $W'_{p,t}$ and $\Phi_{p,t}(W'_{p,t})$
  are both contained in
  the $C^\infty$ manifold $M^* = M \setminus S_{\delta_n}$ without boundary.
\end{lem}

\begin{proof}
  Let us first prove
  $\Phi_{p,t}(W'_{p,t}) \subset M \setminus (\Cut_p \cup S_{\delta_n})$.
  It is clear that $\Phi_{p,t}(W'_{p,t}) \subset M \setminus \Cut_p$.
  To prove $\Phi_{p,t}(W'_{p,t}) \subset M \setminus S_{\delta_n}$,
  we take any point $x \in W'_{p,t}$.
  Since $\Phi_{p,t}(x)$ is not a cut point of $p$ and
  by Lemma \ref{lem:geod}, $\Phi_{p,t}(x)$ is not $\delta_n$-singular.
  Therefore, $\Phi_{p,t}(W'_{p,t}) \subset M \setminus (\Cut_p \cup S_{\delta_n})$.

  Let us next prove $\Phi_{p,t}(W'_{p,t}) \supset
  M \setminus (\Cut_p \cup S_{\delta_n})$.
  Take any point $y \in M \setminus (\Cut_p \cup S_{\delta_n})$ and
  join $p$ to $y$ by a minimal geodesic $\gamma : [\,0,1\,] \to M$.
  Then, $\Phi_{p,t}(\gamma(t)) = y$.  Since $y \not\in \Cut_p$,
  the geodesic $\gamma$ is unique and so $\Phi_{p,t}|_{W'_{p,t}}$ is injective.
  By Lemma \ref{lem:geod}, $\gamma(t) = (\Phi_{p,t}|_{W'_{p,t}})^{-1}(y)$
  is not $\delta_n$-singular and belongs to $W'_{p,t}$.
  This completes the proof.
\end{proof}

By the local Lipschitz continuity of $\Phi_{p,t}$
and by 3.1.8 of \cite{Fd:gmeas},
$\Phi_{p,t}|_{W'_{p,t}}$ is approximately differentiable
$\Hm^n$-a.e.~on $W'_{p,t}$.
The following lemma is essential for the proof of Theorem \ref{thm:BG}.

\begin{lem} \label{lem:Jac}
  Let $p \in M$ and $0 < t < 1$.
  Then, the approximate Jacobian determinant of $\Phi_{p,t}|_{W'_{p,t}}$
  satisfies that
  \[
  |\det\ap d(\Phi_{p,t}|_{W'_{p,t}})_x|
  \le \frac{s_\kappa(r_p(x)/t)^{n-1}}{t\,s_\kappa(r_p(x))^{n-1}}
  \]
  for any approximately differentiable point $x \in W'_{p,t} \setminus S_M$
  of $\Phi_{p,t}|_{W'_{p,t}}$.
\end{lem}

\begin{proof}
  Let $x \in W'_{p,t} \setminus S_M$ be an approximately differentiable
  point of $\Phi_{p,t}|_{W'_{p,t}}$ and let $\epsilon > 0$ be a small number.
  Note that $K_xM$ and $K_{\Phi_{p,t}(x)}M$ are both isometric to $\R^n$
  and identified with the tangent spaces.
  We take two charts $(U,\varphi)$ and $(V,\psi)$ of $M \setminus S_{\delta_n}$
  around $x$ and $\Phi_{p,t}(x)$ respectively such that
  $|\;|\varphi(y)-\varphi(z)|/d(y,z) - 1\,| < \epsilon$ for any different
  $y,z \in U$ and $\psi$ satisfies the same inequality on $V$.
  In particular, every eigenvalue of the differentials
  $d\varphi_x : K_xM \to \R^n$ and
  $d\psi_{\Phi_{p,t}(x)} : K_{\Phi_{p,t}(x)}M \to \R^n$
  is between $1-\epsilon$ and $1+\epsilon$.
  Put
  \begin{align*}
    \bar\Phi &:= \psi\circ\Phi_{p,t}|_{W'_{p,t}}\circ\varphi^{-1} :
    \varphi(W'_{p,t} \cap U) \to \psi(V),\\
    \bar x &:= \varphi(x), \qquad
    L := \ap d\bar\Phi_{\bar x} : \R^n \to \R^n.
  \end{align*}
  For simplicity we set $D := \ap d(\Phi_{p,t}|_{W'_{p,t}})_x
  : K_xM \to K_{\Phi_{p,t}(x)}M$.
  Then,
  \[
  D = (d\psi_{\Phi_{p,t}(x)})^{-1} \circ L \circ d\varphi_x.
  \]
  By the definition of the approximate differential,
  for any $r > 0$ with $B(x,r) \subset U$,
  the set of $\bar y \in B(\bar x,r)$ satisfying
  \[
  |\,\bar\Phi(\bar y) - \bar\Phi(\bar x) - L(\bar y - \bar x)\,|
  \ge \epsilon\, |\,\bar x - \bar y\,|
  \]
  has $\Hm^n$-measure $\le \Theta(r|\bar\Phi,\bar x)\, \Hm^n(B(\bar x,r))$,
  where $B(\bar x,r)$ is a Euclidean metric ball.
  Take any $u \in \Sigma_xM$ and fix it.  Let $r > 0$ be any number.
  We set
  \[
  C(u,r,\epsilon) := \{\; v \in B(o_x,r) \setminus \{o_x\} \subset K_xM
  \mid \angle(u,v) < \epsilon\;\}.
  \]
  It follows from Lemma \ref{lem:HmKp}(1) that
  \begin{align*}
    &\Hm^n(\varphi(\exp_x(C(u,r/2,\epsilon) \cap U_x)))\\
    &\ge (1-\epsilon)^n\,\Hm^n(\exp_x(C(u,r/2,\epsilon) \cap U_x))\\
    &\ge (1-\epsilon)^n\,(\Hm^n(C(u,1/2,\epsilon))
      - \Theta(r|x,n))\,r^n.
  \end{align*}
  Since $\Hm^n(C(u,1/2,\epsilon))$ is positive, we have
  \[
  \lim_{r \to 0}
  \frac{\Hm^n(\varphi(\exp_x(C(u,r/2,\epsilon) \cap U_x)))}
  {\Hm^n(B(\bar x,r))} > 0.
  \]
  Note that $\varphi(\exp_x(C(u,r/2,\epsilon) \cap U_x))$ is contained in
  $B(\bar x,r)$ because $\epsilon$ is small enough.
  Therefore, supposing $r \ll \epsilon$,
  there is a point $\bar y \in B(\bar x,r)$ such that
  \begin{align*}
    &\bar y \in \varphi(\exp_x(C(u,r/2,\epsilon) \cap U_x)),\\
    &|\,\bar\Phi(\bar y) - \bar\Phi(\bar x) - L(\bar y - \bar x)\,|
    < \epsilon\, d(\bar x,\bar y).
  \end{align*}
  Setting $y := \varphi^{-1}(\bar y)$ and $v_{xy} := (\exp_x|_{U_x})^{-1}(y)$,
  we have $\angle(u,v_{xy}) < \epsilon$.
  For simplicity we write 
  $a \le (1+\Theta(\epsilon|p,t,x))\,b + \Theta(\epsilon|p,t,x)$
  by $a \lesssim b$.  Note that since $r \ll \epsilon$,
  all $\Theta(r|\cdots)$ become $\Theta(\epsilon|\cdots)$.
  Since $|v_{xy}| = d(x,y)$ and
  $|d\varphi_x(v_{xy})-(\bar y-\bar x)| \le \Theta(\epsilon|x)\,d(x,y)$
  (cf.~Lemma 3.6(2) of \cite{OS:rstralex}), we have
  \begin{align*}
    |D(u)| &\lesssim |D(v_{xy}/|v_{xy}|)|
    \lesssim
    \frac{|L(\bar y-\bar x)|}{d(x,y)}\\
    &\lesssim \frac{|\bar\Phi(\bar y) - \bar\Phi(\bar x)|}{d(x,y)}
    \lesssim \frac{d(\Phi_{p,t}(x),\Phi_{p,t}(y))}{d(x,y)}.
  \end{align*}

  We are going to estimate the last formula.
  Denote by $M^2(\kappa)$ a complete simply connected $2$-dimensional
  space form of curvature $\kappa$.
  We take three points $\tilde p, \tilde x, \tilde y \in M^2(\kappa)$
  such that $d(\tilde p,\tilde x) = d(p,x)$, $d(\tilde p,\tilde y) = d(p,y)$, 
  and $d(\tilde x,\tilde y) = d(x,y)$.
  The triangle comparison condition tells that
  $d(\Phi_{p,t}(x),\Phi_{p,t}(y)) \le
  d(\Phi_{\tilde p,t}(\tilde x),\Phi_{\tilde p,t}(\tilde y))$,
  where $\Phi_{\tilde p,t}$ is the radial expansion on $M^2(\kappa)$.
  Since $d(\tilde x,\tilde y) = d(x,y) < r \ll \epsilon$, we have
  \begin{align*}
    \frac{d(\Phi_{\tilde p,t}(\tilde x),\Phi_{\tilde p,t}(\tilde y))}
    {d(\tilde x,\tilde y)}
    \lesssim |d(\Phi_{\tilde p,t})_{\tilde x}(v_{\tilde x\tilde y}/|v_{\tilde x\tilde y}|)|.
  \end{align*}
  Let $\tilde\gamma$ be the minimal geodesic from $\tilde p$
  passing through $\tilde x$.
  We denote by $\tilde\theta$ the angle between $v_{\tilde x\tilde y}$ and
  $\tilde\gamma'(t_{\tilde x})$,
  where $t_{\tilde x}$ is taken in such a way that
  $\tilde\gamma(t_{\tilde x}) = \tilde x$.
  Set
  \[
  \lambda(\xi) := \sqrt{\frac{1}{t^2}\cos^2\xi
    + \frac{s_\kappa(r_p(x)/t)^2}{s_\kappa(r_p(x))^2}\sin^2\xi},
  \qquad \xi \in \R.
  \]
  A calculation using Jacobi fields yields
  $|d(\Phi_{\tilde p,t})_{\tilde x}(v_{\tilde x\tilde y}/|v_{\tilde x\tilde y}|)|
  = \lambda(\tilde\theta)$.
  Combining the above estimates, we have
  \[
  |D(u)| \lesssim \lambda(\tilde\theta).
  \]
  Let $\gamma$ be the minimal geodesic from $p$ passing through $x$
  and let $t_x$ be a number such that $\gamma(t_x) = x$.
  Denote by $\theta$ the angle between $v_{xy}$ and $\gamma'(t_x)$ and
  by $\theta_u$ the angle between $u$ and $\gamma'(t_x)$.
  It follows from $\angle(u,v_{xy}) < \epsilon$
  that $|\theta - \theta_u| < \epsilon$.
  By 5.6 of \cite{BGP} we have
  $|\,\theta - \tilde\theta\,| \le \Theta(r|p,t,x) \le \Theta(\epsilon|p,t,x)$.
  Therefore we have
  $|D(u)| \lesssim \lambda(\theta_u)$.
  Taking the limit as $\epsilon \to 0$ yields that
  \[
  |D(u)| \le \lambda(\theta_u)
  \]
  for any $u \in \Sigma_xM$, which together with Hadamard's inequality implies
  \[
  |\det D| \le \lambda(0)\,\lambda(\pi/2)^{n-1}
  = \frac{s_\kappa(r_p(x)/t)^{n-1}}{t\,s_\kappa(r_p(x))^{n-1}}.
  \]
  This completes the proof of Lemma \ref{lem:Jac}.
\end{proof}

\begin{proof}[Proof of Theorem \ref{thm:BG}]
  For the proof, it suffices to prove that
  \begin{align}
    \int_{W_{p,t}} f\circ\Phi_{p,t}(x)\, d\Hm^n(x) 
    \ge \int_M f(y)\,\frac{t\, s_\kappa(t\,r_p(y))^{n-1}}{s_\kappa(r_p(y))^{n-1}}
    \,d\Hm^n(y)
    \label{eq:BG}
  \end{align}
  for any $\Hm^n$-measurable function $f : M \to [\,0,+\infty\,)$
  with compact support.
  Since $\Phi_{p,t}|_{W'_{p,t}} : W'_{p,t} \to M \setminus (\Cut_p \cup S_{\delta_n})$
  is bijective,
  the area formula (cf.~3.2.20 of \cite{Fd:gmeas}) implies that
  \begin{align} \label{eq:area}
    &\int_{W'_{p,t}} F\circ\Phi_{p,t}(x)\,
    |\det\ap d(\Phi_{p,t}|_{W'_{p,t}})_x|\;d\Hm^n(x)\\
    &= \int_{M \setminus (\Cut_p \cup S_{\delta_n})} F(y)\;d\Hm^n(y) \notag
  \end{align}
  for any $\Hm^n$-measurable function $F : M \to [\,0,+\infty\,)$
  with compact support.
  We set
  \[
  F(y) := f(y)\,\frac{t\, s_\kappa(t\,r_p(y))^{n-1}}{s_\kappa(r_p(y))^{n-1}},
  \quad y \in M \setminus \Cut_p,
  \]
  in \eqref{eq:area}.  Then, since $\Hm^n(\Cut_p) =
  \Hm^n(S_{\delta_n}) = 0$ and by Lemma \ref{lem:Jac}, we obtain \eqref{eq:BG}.
  This completes the proof of the theorem.
\end{proof}

\bibliographystyle{amsplain}
\bibliography{iwpt}

\begin{thebibliography}{10}

\bibitem{BGP}
Yu. Burago, M.~Gromov, and G.~Perel'man, \emph{A. {D}. {A}leksandrov spaces
  with curvatures bounded below}, Uspekhi Mat. Nauk \textbf{47} (1992),
  no.~2(284), 3--51, 222, translation in Russian Math. Surveys 47 (1992), no.
  2, 1--58.

\bibitem{CC:strRicI}
J.~Cheeger and T.~H. Colding, \emph{On the structure of spaces with {R}icci
  curvature bounded below. {I}}, J. Differential Geom. \textbf{46} (1997),
  no.~3, 406--480.

\bibitem{Fd:gmeas}
H.~Federer, \emph{Geometric measure theory}, Springer, Berlin, 1969.

\bibitem{KMS:lap}
K.~Kuwae, Y.~Machigashira, and T.~Shioya, \emph{Sobolev spaces, {L}aplacian,
  and heat kernel on {A}lexandrov spaces}, Math. Z. \textbf{238} (2001), no.~2,
  269--316.

\bibitem{KwSy:splitting}
K.~Kuwae and T.~Shioya, \emph{A topological splitting theorem for weighted
  {A}lexandrov spaces}, preprint.

\bibitem{KwSy:geneMCP}
\bysame, \emph{On generalized measure contraction property and energy
  functionals over {L}ipschitz maps}, Potential Anal. \textbf{15} (2001),
  no.~1-2, 105--121, ICPA98 (Hammamet).

\bibitem{KwSy:sobmet}
\bysame, \emph{Sobolev and {D}irichlet spaces over maps between metric spaces},
  J. Reine Angew. Math. \textbf{555} (2003), 39--75.

\bibitem{KwSy:lapcomp}
\bysame, \emph{{L}aplacian comparison for {A}lexandrov spaces}, preprint, 2007.

\bibitem{LV:Ricmm}
J.~Lott and C.~Villani, \emph{Ricci curvature for metric-measure spaces via
  optimal transport}, to appear in Ann. Math., 2006.

\bibitem{Ota:mcp}
S.~Ohta, \emph{On the measure contraction property of metric measure spaces},
  Comment. Math. Helv. \textbf{82} (2007), no.~4, 805--828.

\bibitem{Ot:secdiff}
Y.~Otsu, \emph{Almost everywhere existance of second differentiable structure
  of {A}lexandrov spaces}, preprint.

\bibitem{OS:rstralex}
Y.~Otsu and T.~Shioya, \emph{The {R}iemannian structure of {A}lexandrov
  spaces}, J. Differential Geom. \textbf{39} (1994), no.~3, 629--658.

\bibitem{Pr:DC}
G.~Perelman, \emph{{DC}-structure on {A}lexandrov space}, preprint.

\bibitem{Pt:parallel}
A.~Petrunin, \emph{Parallel transportation for {A}lexandrov space with
  curvature bounded below}, Geom. Funct. Anal. \textbf{8} (1998), no.~1,
  123--148.

\bibitem{Rm:Poincare}
A.~Ranjbar-Motlagh, \emph{Poincar\'e inequality for abstract spaces}, Bull.
  Austral. Math. Soc. \textbf{71} (2005), no.~2, 193--204.

\bibitem{Sy:raysalex}
T.~Shioya, \emph{Mass of rays in {A}lexandrov spaces of nonnegative curvature},
  Comment. Math. Helv. \textbf{69} (1994), no.~2, 208--228.

\bibitem{Sy:Ronsetsu}
\bysame, \emph{Geometric analysis on Alexandrov spaces}, 
to appear in Sugaku Expositions. 

\bibitem{St:heat}
K.-T. Sturm, \emph{Diffusion processes and heat kernels on metric spaces}, Ann.
  Probab. \textbf{26} (1998), no.~1, 1--55.

\bibitem{St:geomI}
\bysame, \emph{On the geometry of metric measure spaces. {I}}, Acta Math.
  \textbf{196} (2006), no.~1, 65--131.

\bibitem{St:geomII}
\bysame, \emph{On the geometry of metric measure spaces. {II}}, Acta Math.
  \textbf{196} (2006), no.~1, 133--177.

\bibitem{Wt:localcut}
M.~Watanabe, \emph{Local cut points and metric measure spaces with {R}icci
  curvature bounded below}, Pacific J. Math. \textbf{233} (2007), no.~1,
  229--256.

\end{thebibliography}

\end{document}